\documentclass[10pt,a4paper]{article}

\usepackage[T1]{fontenc}
\usepackage{amsmath}
\usepackage{amsfonts}
\usepackage{amssymb}
\usepackage{amsthm}
\usepackage[headsep=10mm]{geometry}
\usepackage{caption}
\usepackage{subcaption}
\usepackage[pdftex]{graphicx}
\usepackage{bm}
\usepackage[hidelinks]{hyperref}
\usepackage{fancyhdr}
\usepackage{dsfont}
\usepackage{stmaryrd}
\usepackage{enumitem}
\usepackage{color}
\usepackage{authblk}

\numberwithin{equation}{section}

\newcommand{\pfrac}[2]{\frac{\partial{#1}}{\partial{#2}}}
\newcommand{\dx}[1]{\hspace{0.75mm}\mathrm{d}{#1}}

\newcommand{\curlybrackets}[1]{\left\lbrace{#1}\right\rbrace}

\definecolor{purple}{rgb}{1,0,1}

\begin{document}
\author[1,2]{Thomas M. Bendall}
\author[1]{Colin J. Cotter}
\author[1]{Darryl D. Holm}
\affil[1]{Mathematics Department, Imperial College, London, UK}
\affil[2]{Dynamics Research, Met Office, Exeter, UK}
\title{Perspectives on the Formation of Peakons in the Stochastic Camassa-Holm Equation}
\date{}
\maketitle

\begin{abstract}
\noindent
A famous feature of the Camassa-Holm equation is its admission of peaked soliton solutions known as peakons.
We investigate this equation under the influence of stochastic transport.
Noting that peakons are weak solutions of the equation, we present a finite element discretisation for it, which we use to explore the formation of peakons. \\
\\
Our simulations using this discretisation reveal that peakons can still form in the presence of stochastic perturbations.
Peakons can emerge both through wave breaking, as the slope turns vertical, and without wave breaking as the inflection points of the velocity profile rise to reach the summit.
\end{abstract}
\section{Introduction}
\subsection{Overview}
We are dealing with stochastic perturbations of partial differential equations (PDEs) for fluid dynamics, whose deterministic solutions develop singularities by producing a vertical derivative, which is commonly known as \emph{wave breaking}.
A famous example is the Burgers equation, which has been studied with a stochastic component several times, for instance by \cite{bertini1994stochastic} and \cite{alonso2019burgers}. In the latter case, the stochastic perturbations take the form of multiplicative noise, which will be our focus here.
A striking feature is that the addition of stochastic perturbations has been found to affect the  development of singularities; examples include \cite{flandoli2010well} and \cite{flandoli2011random}.\\
\\
Here we consider the effect of stochastic perturbations on a well-known equation that can also exhibit wave breaking and singularity formation:  the Camassa-Holm (CH) equation of \cite{camassa1993integrable}, which has recently been studied in the stochastic case by \cite{holm2016variational} and \cite{crisan2018wave}.
This stochastic equation can be expressed as
\begin{equation}
\mathrm{d}m(x,t) + (\partial_x m + m \partial_x)(u\dx{t} + \varXi(x)\circ \mathrm{d}W) = 0
\,,\quad\mathrm{with}\quad 
m = u - \alpha^2u_{xx},
\label{eqn:SCH-eqn}
\end{equation}
for the evolution in time $t\in\mathbb{R}$ of the fluid momentum density $m(x,t)$ and the velocity $u(x,t)$ on the real line $x\in\mathbb{R}$.
The constant $\alpha$ has dimensions of length.
In equation (\ref{eqn:SCH-eqn}), the time differential $\mathrm{d}$ is short notation for a stochastic integral; the spatial partial derivative $\partial_x$ acts rightwards on all products; the symbol $\circ$ denotes a Stratonovich stochastic process; $\mathrm{d}W$ is a Brownian motion with spatially modulated amplitude $\varXi(x)$; and one interprets the effect of the stochasticity as a noisy perturbation of the transport velocity.
See \cite{flandoli2011random} for more discussion of such stochastic transport.
In the deterministic case, with $\varXi(x)=0$, it is well-known that (\ref{eqn:SCH-eqn}) admits weak singular solutions known as \emph{peakons}, which refers to their shape with a vertical derivative at their peak. \\
\\
The paper begins by presenting a finite element discretisation to a form of (\ref{eqn:SCH-eqn}), inspired by the acknowledgement that peakons are weak solutions.
We then use this discretisation to study the formation of peakons in deterministic and stochastic cases.
We look at the formation of peakons both with and without wave breaking.
In this second case, rather than the slope turning vertical the inflection points of the profile rise until they reach the summit.
With the addition of stochastic noise, our simulations suggest that peakon formation (with and without wave breaking) persists.
\subsection{A stochastic variational principle}\label{sec:stochastic variational principle}
Stochastic perturbations are of current interest in the modelling of geophysical fluids, where they can express the uncertainty in the effects of the unresolved processes upon the resolved flow.
In this context, a new class of stochastic fluid equations was introduced by \cite{holm2015variational}, which crucially preserves many of the circulation properties of the respective deterministic equations.
For instance in the stochastic quasi-geostrophic equations studied by \cite{bendall2018statistical} and \cite{cotter2019numerically}, the potential vorticity is still preserved by stochastic material transport.
These equations are derived from a stochastically constrained variational principle $\delta S=0$, where the action integral $S$ is given by
\begin{equation}
S = (\bm{u}, p, q) = \int \left[  \ell(\bm{u}, q) \dx{t}  + 
\left\langle p, \mathrm{d}q + \mathfrak{L}_{\mathrm{d}\bm{x}_t} q\right \rangle \right],
\end{equation}
where $\bm{u}$ is the velocity advecting quantity $q$ and $\ell(\bm{u},q)$ is the Lagrangian of the deterministic fluid.
The Lagrange multiplier $p$ enforces the stochastic transport of $q$ by the Lie derivative $\mathfrak{L}_{\mathrm{d}\bm{x}_t}$ along a vector field $\mathrm{d}\bm{x}_t$, which is given by
\begin{equation}
\mathrm{d}\bm{x}_t = \bm{u}(\bm{x},t) \dx{t} - \sum_j \bm{\varXi}^j(\bm{x})\circ \mathrm{d}W^j(t), \label{eqn:dxt}
\end{equation}
with the sum over multiple Wiener processes, which are modulated in space by a series of functions $\bm{\varXi}^j(\bm{x})$, indexed by superscript $j$.
The angled brackets denote the spatial integral over the domain $\varOmega$ of the pairing of $p$ and $q$,
\begin{equation}
\left\langle p, q \right\rangle = \int_\varOmega p q \dx{x}.
\end{equation}
A thorough description of this formulation can be found in \cite{holm2015variational}.
\subsection{The deterministic and stochastic Camassa-Holm equations}
The deterministic Camassa-Holm (CH) equation of \cite{camassa1993integrable} can be expressed as
\begin{equation}\label{eqn:ch hydrodynamic}
u_t =
-\partial_x\curlybrackets{\frac{1}{2}u^2+ K\ast \left[u^2+\frac{\alpha^2}{2}\left(u_x\right)^2\right]},
\end{equation}
which is sometimes known as the \emph{advective} or \emph{hydrodynamic} form of the equation.
In contrast to (\ref{eqn:SCH-eqn}), this does not contain double spatial derivatives of $u$.
For more discussion of this form, see for example \cite{constantin1998global} and \cite{crisan2018wave}.
The operator $K$ is the Green's function for the Helmholtz operator $\left(1-\alpha^2\partial_{xx}\right)$ that relates the momentum density $m$ and the velocity $u$, so that
\begin{equation} \label{eqn:greens function}
K\ast m = K\ast\left(u-\alpha^2 u_{xx}\right) = u.
\end{equation}
The CH equation (\ref{eqn:ch hydrodynamic}) was derived at one order beyond the celebrated Korteweg-de Vries (KdV) equation in the asymptotic expansion of the Euler fluid equations for non-linear shallow water waves on a free surface propagating under the restoring force of gravity \cite{dullin2001integrable, dullin2003camassa,dullin2004asymptotically}.
It is integrable \cite{camassa1993integrable, constantin2001scattering} and has a bi-Hamiltonian structure \cite{fuchssteiner1981symplectic,fuchssteiner1996some}.\\
\\
The stochastic Camassa-Holm (SCH) equation (\ref{eqn:SCH-eqn}) was derived by \cite{holm2016variational}  from the stochastic variational principle of \cite{holm2015variational}  discussed in Section \ref{sec:stochastic variational principle}.
As seen in \cite{crisan2018wave}, (\ref{eqn:SCH-eqn}) also possesses a hydrodynamic form:
\begin{equation}\label{eqn:sch hydro}
\begin{split}
\mathrm{d}u = & -\partial_x\curlybrackets{\frac{1}{2}u^2 + K\ast \left[u^2+\frac{\alpha^2}{2}\left(u_x\right)^2\right]}\mathrm{d}t \\
& - \sum_j \curlybrackets{u_x\varXi^j+ K\ast \left[2u\varXi_x^j+\alpha^2u_x\varXi^j_{xx}\right]}\circ \mathrm{d}W^j.
\end{split}
\end{equation}
Unlike (\ref{eqn:ch hydrodynamic}), the right hand side of (\ref{eqn:sch hydro}) cannot be written purely as a gradient, but it also contains no double derivatives of $u$.
Other stochastic Camassa-Holm equations have also been studied, including that of \cite{chen2012well} which used additive noise, and that of \cite{lv2020dependence} which used a type of multiplicative noise which differs from that derived from the variational principle of \cite{holm2015variational}.
Equation (\ref{eqn:sch hydro}) will be the focus for the rest of this paper.\\
\\
The deterministic CH equation (\ref{eqn:ch hydrodynamic}) admits solutions with singularities  known as \emph{peakons}, which possess a sharp peak at the apex of their velocity profile.
In fact, due to the singularity in the spatial derivative of the peakon solutions, they should be interpreted as \textit{weak} solutions to the CH equation (\ref{eqn:ch hydrodynamic}), as discussed by \cite{constantin1998global,constantin1998wave,constantin2000global}.
The shape of the peakon profile turns out to be the Green's function $K$.
This can be explained by recalling that the peakon solutions for the momentum density in the deterministic form of (\ref{eqn:SCH-eqn}) may be expressed by the singular momentum map \cite{holm2005momentum},
\begin{align}
m(x,t) = \sum_k p_k(t)\delta\left(\frac{x-q_k(t)}{\alpha}\right)
\quad\mathrm{so}\quad
u(x,t) = \sum_k p_k(t)K\ast\delta\left(\frac{x-q_k(t)}{\alpha}\right)
\,,
\label{eqn:sing-momap}
\end{align}
in which $K(x-q_k(t)) = \frac{1}{2}\exp(-|x-q_k(t)|/\alpha)$ is the Green's function on the real line with homogeneous boundary conditions.
Remarkably, the purely discrete spectrum of the isospectral eigenvalue equation, which demonstrates that the CH equation is a completely integrable Hamiltonian soliton system, also implies that \emph{only} the singular solutions in (\ref{eqn:sing-momap}) will persist.
The continuation of peakon solutions beyond their emergence has been investigated by studies including \cite{bressan2007global} and \cite{holden2007global}.
For such a sum of peakons, the canonical coordinates giving the peakon momentum $p_k$ and position $q_k$ obey the following coupled ordinary differential equations:
\begin{equation}
\frac{\mathrm{d}p_k}{\mathrm{d}t}=-p_k\pfrac{u(q_k)}{q_k} \quad \mathrm{and} \quad \frac{\mathrm{d}q_k}{\mathrm{d}t}=u(q_k).
\end{equation}
The stochastic equation (\ref{eqn:sch hydro}) also admits weak peakon solutions via (\ref{eqn:sing-momap}), as mentioned by \cite{holm2016variational}.
\cite{holm2016variational} derived the stochastic (ordinary) differential equations obeyed by these peakons:
\begin{equation}\label{eqn:peakon differentials raw}
\mathrm{d}p_k =- p_k(t)\pfrac{\mathrm{d}x_t(q_k(t))}{q_k}  \quad
\mathrm{and} \quad \mathrm{d}q_k = \mathrm{d}x_t(q_k(t)),
\end{equation}
where
\begin{equation}
\mathrm{d}x_t(q_k(t)) = u(q_k(t))\dx{t} + \sum_j \varXi^j(q_k(t))\circ \mathrm{d}W^j.
\end{equation}
Evaluating (\ref{eqn:peakon differentials raw}) gives
\begin{subequations} \label{eqn:diff peakons Strat non-periodic}
\begin{align} 
\mathrm{d}p_k & = -\frac{p_k}{2\alpha} \sum_l p_l \ \mathrm{sign}(q_k-q_l) e^{-|q_k-q_l|/\alpha}\dx{t} - p_k\sum_j\varXi^j_x(q_k)\circ \mathrm{d}W^j, \\
\mathrm{d}q_k & = \frac{1}{2}\sum_l p_l \ e^{-|q_k-q_l|/\alpha}\dx{t} + \sum_j \varXi^j(q_k)\circ \mathrm{d}W^j.
\end{align}
\end{subequations}
where we take $\mathrm{sign}(0)=0$.
\cite{holm2016variational} discretised (\ref{eqn:diff peakons Strat non-periodic}) using a variational integrator and used it to investigate the interaction of peakons.
This was followed by \cite{holm2018stochastic} and \cite{holm2018new} which further explored stochastic Hamiltonian systems.
This paper continues in the same spirit as \cite{holm2016variational} of using numerical methods to understand the equation, as we present a finite element discretisation for (\ref{eqn:sch hydro}) and using it to investigate the formation of peakons in the SCH equation.
\\
\\
A final property to mention concerning peakons is the \emph{steepening lemma} of \cite{camassa1993integrable}, which is repeated for instance in \cite{crisan2018wave}.
This demonstrates that peakons emerge from a class of smooth, confined initial velocity profiles.
Such profiles have an inflection point with a negative slope $s$ at $x=\widetilde{x}$ to the right of their maxima.
The steepening lemma of \cite{camassa1993integrable}, as well as the argument of \cite{constantin1998wave}, shows that if $s$ is sufficiently negative then the slope will become vertical within finite time, and thus a peakon emerges through wave breaking.
Considering a velocity profile $u$ with momentum $m=u-\alpha^2u_{xx}$, \cite{mckean1998breakdown} argues that wave breaking will occur if there is some negative $m$ at points $x_{-}$ that lie to the right of points $x_{+}$ at which $m$ is positive.
In a periodic domain this must be the case, unless $m>0$ everywhere.
\\
\\
\cite{crisan2018wave} then investigated the analogous steepening lemma for the SCH equation,  again looking at only a class of initially smooth, confined solutions whose negative slope $s$ at $\widetilde{x}$ was sufficiently negative.
Focusing on the case of a single stochastic basis function, $\varXi(x)=\xi$ (a constant), the authors found that the expectation of the slope at the inflection point also blows up in finite time.
Despite this, for individual realisations of the noise it was unclear whether wave breaking would  always occur, as \cite{crisan2018wave} found that the slope becomes vertical with a positive probability which may not necessarily be unity.
This was interpreted as the probability of peakon formation. \\
\\ 
As we discuss in Section \ref{sec:results}, peakons can still form from smooth, confined initial conditions, even if $s$ is not sufficiently negative for wave breaking to occur.
Instead, the inflection point rises up the velocity profile,
with the peakon emerging as the inflection point reaches the profile's maximum.
This emphasises that the stochastic steepening lemma of \cite{crisan2018wave} describes the probability of wave breaking (for initial profiles with sufficiently negative $s$) and not the probability of peakon formation.\\
\\
Some outstanding questions therefore remain.
While the introduction of stochastic transport into the SCH equation cannot be expected to preserve the complete integrability of the unperturbed CH equation, one may ask whether peakons still form from smooth, confined initial conditions.
We look first at the case of $\varXi(x)=\xi$ (a constant) considered by \cite{holm2016variational} and \cite{crisan2018wave}.
In the situation of \cite{crisan2018wave}, with sufficiently negative $s$, is the probability of wave breaking actually unity or merely non-zero?
What about the probability of peakon formation?
When the initial condition has a shallower slope, do peakons still form under stochastic perturbation, and can wave breaking now occur even when it wouldn't in the deterministic case?
The purpose of this paper is investigate these questions numerically via a finite element discretisation for (\ref{eqn:sch hydro}).
\\
\\
Here is the paper's structure. In Section \ref{sec:weak peakons}, we verify that peakons do indeed satisfy the stochastic Camassa-Holm equation by writing it in hydrodynamic form.
In Section \ref{sec:discretisation} we present a finite element discretisation for the stochastic Camassa-Holm equation, showing that it numerically converges to the stochastic (ordinary) differential equations for peakons in Section \ref{sec:numerical convergence}.
In Section \ref{sec:results}, we numerically investigate the formation of peakons using our discretisation.

%%%%%%%%%%%%%%%%%%%%%%%%%%%%%%%%%%%%%%%%%%%
\section{Peakon Solutions to the Stochastic Camassa-Holm equation} \label{sec:weak peakons}
In this section we verify that peakons are also solutions to the stochastic Camassa-Holm equation.
These peakons take the form
\begin{equation} \label{eqn:u in terms of K}
u = \sum_k p_k(t) K\ast\delta\left(\frac{x-q_k(t)}{\alpha}\right),
\end{equation}
in which $q_k(t)$ represents the position of the peakons and $p_k(t)$ represents their momentum, and together they form pairs of canonical coordinates satisfying a Hamiltonian system.
Although this was described in \cite{crisan2018wave}, here we emphasise that peakons are \textit{weak} solutions to (\ref{eqn:sch hydro}), since the derivative $u_x$ is not defined at the peak.
In other words, peakons are solutions to
\begin{equation} \label{eqn:weak sch-hydro}
\begin{split}
\int_\varOmega\phi \ \mathrm{d}u \dx{x} = & \int_\varOmega \phi_x \curlybrackets{\frac{1}{2}u^2 + K\ast \left[u^2+\frac{\alpha^2}{2}\left(u_x\right)^2\right]}\mathrm{d}t \dx{x} \\
& - \sum_j \int_\varOmega \phi \curlybrackets{u_x\varXi^j + K\ast \left[2u\varXi_x^j+\alpha^2u_x\varXi^j_{xx}\right]}\circ \mathrm{d}W^j \ \mathrm{d}x,
\end{split}
\end{equation}
for all $\phi \in H^1(\varOmega)$.\\
\\
In later sections we numerically solve the stochastic Camassa-Holm equation in the periodic interval $x\in [0,L]$ so we take this here for our domain $\varOmega$.
In this case, the operator $K(x)$ is given by
\begin{equation} \label{eqn:greens function}
K(x) = \frac{1}{2\left(1-e^{-L/\alpha}\right)}\left[e^{-\mathrm{mod}(x,L)/\alpha}+e^{-L/\alpha}e^{\mathrm{mod}(x,L)/\alpha} \right].
\end{equation}
On the periodic interval, the velocity for series of peakons with momenta $p_k$ centred at positions $q_k$ is given by
\begin{equation} \label{eqn:periodic peakon}
u = \sum_k u^k =\frac{1}{2}\sum_k\frac{p_k(t)}{1-e^{-L/\alpha}}\left[e^{-\mathrm{mod}(x-q_k(t),L)/\alpha} + e^{-L/\alpha}e^{\mathrm{mod}(x-q_k(t),L)/\alpha} \right],
\end{equation}
where we define $u^k$ to be the contribution from a single peakon.
Other works using periodised peakons include \cite{constantin1999shallow}, \cite{xu2008local} and \cite{holden2008numerical}.
This reduces to the usual set of peakons on $\mathbb{R}$ when $L\to\infty$:
\begin{equation} \label{eqn:old peakon}
u = \frac{1}{2}\sum_k p_k(t)e^{-|x-q_k(t)|/ \alpha}.
\end{equation}
In the periodic case, the stochastic ODEs for the evolution of the momenta $p_k$ and positions $q_k$ are given by
\begin{subequations}\label{eqn:diff peakons Strat}
\begin{align} 
\begin{split}
\mathrm{d}p_k & = \\
& -\frac{p_k}{2\alpha(1-e^{-L/\alpha})} \sum_l p_l \ \mathrm{sign}(q_k-q_l) \left[e^{-\mathrm{mod}(q_k-q_l,L)/\alpha} - e^{-L/\alpha}e^{\mathrm{mod}(q_k-q_l),L)/\alpha}\right]\dx{t} \\
& - p_k\sum_j\varXi^j_x(q_k)\circ \mathrm{d}W^j,
\end{split} \\
\begin{split}
\mathrm{d}q_k & = \frac{1}{2(1-e^{-L/\alpha})}\sum_l p_l \ \left[e^{-\mathrm{mod}(q_k-q_l,L)/\alpha}+e^{-L/\alpha}e^{\mathrm{mod}(q_k-q_l,L)/\alpha}\right]\dx{t} \\
& + \sum_j \varXi^j(q_k)\circ \mathrm{d}W^j.
\end{split}
\end{align}
\end{subequations}
which is the periodic equivalent of (\ref{eqn:diff peakons Strat non-periodic}).
Taking the differential of (\ref{eqn:periodic peakon}) and substituting in equation (\ref{eqn:peakon differentials raw})
yields
\begin{subequations}
\begin{align}
\mathrm{d}u&= \sum_k\left[\pfrac{u^k}{p_k}\dx{p_k} + \pfrac{u^k}{q_k} \dx{q_k}\right], \\
& = -\sum_k\left[\left(u^k \pfrac{u(q_k)}{q_k}+u^k_x u(q_k) \right)\dx{t} + \sum_j \left(u^k\varXi_x^j(q_k)+u^k_x\varXi^j(q_k) \right)\circ\mathrm{d}W^j\right]. \label{eqn:du}
\end{align}
\end{subequations}
In order to verify that (\ref{eqn:periodic peakon}) satisfies (\ref{eqn:weak sch-hydro}), we will substitute it and (\ref{eqn:du}) into (\ref{eqn:weak sch-hydro}) and show that the deterministic and stochastic components of the equation will both vanish separately.
First, substituting just (\ref{eqn:du}) into the left-hand side of (\ref{eqn:weak sch-hydro}) and grouping terms, with the deterministic parts on the left-hand side and stochastic terms on the right-hand side, gives:
\begin{equation} \label{eqn:peakon det and stoch}
\begin{split}
&- \int_\varOmega \phi \sum_k \left(u^k \pfrac{u(q_k)}{q_k}+u^k_x u(q_k) \right) \dx{t} \dx{x} - \int_\varOmega \phi_x \curlybrackets{\frac{1}{2}u^2 + K\ast \left[u^2+\frac{\alpha^2}{2}\left(u_x\right)^2\right]}\mathrm{d}t \dx{x} = \\
& \int_\varOmega \phi \sum_j \curlybrackets{\sum_k \left[u^k\varXi^j_x(q_k) + u^k_x\varXi^j(q_k)\right] - u_x \varXi^j - K\ast\left(2u\varXi_x^j + \alpha^2u_x\varXi_{xx}^j \right) } \circ \mathrm{d}W^j \dx{x}
\end{split}
\end{equation}
Now we consider the multi-peakon solution for $u$ from (\ref{eqn:periodic peakon}).
Although we don't give explicit proof here (for more details see \cite{constantin2000global} and \cite{constantin1999shallow} for the periodic case), this peakon is a weak solution to the deterministic equation, and therefore the left-hand side of (\ref{eqn:peakon det and stoch}) is zero.
Then to show (\ref{eqn:periodic peakon}) is a solution, we must demonstrate that the terms inside the curly brackets $\curlybrackets{\cdot}$ in the stochastic part of (\ref{eqn:peakon det and stoch}) cancel almost everywhere, so that the integral on the right-hand side vanishes. \\
\\
Since $u$ is a linear sum of $u^k$ for each term on the right-hand side of (\ref{eqn:peakon det and stoch}), we only need to consider a single peakon solution $u^k$ and a single $\varXi^j$, taking
\begin{equation}
\varXi^j = \xi\left[C\cos\left(\frac{2\pi j x}{L}\right)+D\sin\left(\frac{2\pi j x}{L}\right) \right],
\end{equation}
as in a periodic domain, any choice of $\varXi$ will be a linear combination of such functions.
With this choice, careful computation of the convolution term in the stochastic part of (\ref{eqn:peakon det and stoch}) reveals that for $x\neq q_k$:
\begin{align} \label{eqn:peakon proof}
\begin{split}
& K\ast\left(2u^k\varXi^j + \alpha^2u^k_x\varXi_{xx}^j \right) =  \\
& \frac{1}{2\alpha(1-e^{-L/\alpha})} \int_0^L \left[e^{-\mathrm{mod}(x-y,L)/\alpha} + e^{-L/\alpha}e^{\mathrm{mod}(x-y,L)/\alpha} \right]\left[ 2u^k(y)\varXi_y^j(y) + \alpha^2u^k_y(y) \varXi_{yy}^j(y) \right] \dx{y} = \\
\begin{split}
& \frac{\xi p_k}{2(1-e^{-L/\alpha})} \left\lbrace
\frac{2\pi j}{L}
\left(e^{-\mathrm{mod}(x-q_k,L) /\alpha} + e^{-L/\alpha}e^{\mathrm{mod}(x-q_k,L)/\alpha} \right)
 \left[-C\sin\left(\frac{2\pi j x}{L}\right)+D\cos\left(\frac{2\pi j x}{L}\right) \right] \right. - \\ 
& \ \ \ \ \ 
 \left. \frac{\mathrm{sign}(x-q_k)}{\alpha}
\left(e^{-\mathrm{mod}(x-q_k,L) /\alpha} - e^{-L/\alpha}e^{\mathrm{mod}(x-q_k,L)/\alpha}  \right)\left[C\cos\left(\frac{2\pi j x}{L}\right)+D\sin\left(\frac{2\pi j x}{L}\right) \right]\right\rbrace
\end{split} \\
& = u^k \varXi_x^j(q_k) + u^k_x\varXi^j(q_k) - u^k_x\varXi^j(x).
\end{split}
\end{align}
Consequently the term inside the curly brackets of the stochastic of (\ref{eqn:peakon det and stoch}) vanishes almost everywhere, and the integral on the right-hand side of (\ref{eqn:peakon det and stoch}) is zero.
Since both the deterministic and stochastic parts of (\ref{eqn:peakon det and stoch}) vanish when the velocity profile is given by (\ref{eqn:periodic peakon}), this periodic peakon is indeed a weak solution of the SCH equation (\ref{eqn:weak sch-hydro}), with the momenta $p_k$ and $q_k$ obeying (\ref{eqn:diff peakons Strat}).\\
\\
Alternatively, taking a constant $\varXi^j=\xi$ in (\ref{eqn:peakon det and stoch}) then the convolution term vanishes as $\varXi^j_x=0$ and $\varXi^j_{xx}=0$.
Since for constant $\xi$, $\varXi^j(x)=\varXi^j(q_k)$, the remaining terms on the right-hand side of (\ref{eqn:diff peakons Strat}) cancel and again we see that (\ref{eqn:periodic peakon}) is a solution of (\ref{eqn:weak sch-hydro}).

\section{A Finite Element Discretisation for the Stochastic Camassa-Holm Equation} \label{sec:discretisation}
As our aim is to numerically investigate the formation of peakons within the stochastic Camassa-Holm equation, we look for weak solutions of the hydrodynamic form (\ref{eqn:sch hydro}), acknowledging that peakons are weak solutions of this equation.
Numerical solutions will then be valid before and after the moment of peakon formation.
While looking for weak solutions, it is natural to discretise the equation using a finite element method.\\
\\
The Camassa-Holm equation has been studied numerically using a range of methods, including finite difference, e.g. \cite{holden2006convergence}, finite volume, e.g. \cite{artebrant2006numerical} and particle methods, e.g. \cite{chertock2012convergence}.
Previous examples of finite element methods to discretise the Camassa-Holm equation include discontinuous Galerkin (DG) methods, such as those of \cite{xu2008local}, \cite{liu2016invariant} and \cite{li2014high}, and the continuous methods of \cite{matsuo2010hamiltonian} and \cite{antonopoulos2019error}.
In contrast to those examples, we discretise the Camassa-Holm equation in its hydrodynamic form. \\
\\
We now present a new finite element discretisation of the stochastic Camassa-Holm equation that we will use in the remainder of the paper.
\\
\\
To obtain the discretisation, we introduce the variables $F$ and $G^j$ and write (\ref{eqn:sch hydro}) as
\begin{subequations}\label{eqn:fem-1}
\begin{align}
\mathrm{d}u & = -\partial_x\curlybrackets{\frac{1}{2}u^2 + F} \dx{t} - \sum_j \curlybrackets{u_x\varXi^j + G^j}\circ \mathrm{d}W^j, \\
F & = K \ast \left[u^2 + \frac{\alpha^2}{2}(u_x)^2\right], \\
G^j &= K \ast \left[2u\varXi^j_x + \alpha^2u_x\varXi^j_{xx} \right].
\end{align}
\end{subequations}
Using that $K$ is the Green's function for the Helmholtz operator, this becomes
\begin{subequations} \label{eqn:fem-2}
\begin{align}
\mathrm{d}u & = -\partial_x\curlybrackets{\frac{1}{2}u^2 + F} \dx{t} - \sum_j \curlybrackets{u_x\varXi^j + G^j}\circ \mathrm{d}W^j, \\
F-\alpha^2F_{xx} & = u^2 + \frac{\alpha^2}{2}(u_x)^2, \\
G^j-\alpha^2G^j_{xx} &= 2u\varXi^j_x + \alpha^2u_x\varXi^j_{xx}.
\end{align}
\end{subequations}
We restrict ourselves to considering a periodic one-dimensional domain $\varOmega$ of length $L_d$.
Now we introduce our finite element space $V$, which is the space of continuous piecewise-linear functions over some partition of $\varOmega$ (often this is known as the continuous Galerkin space CG$_1$).
We multiply each of the equations in (\ref{eqn:fem-2}) by the test functions $\psi$, $\phi$ and $\zeta$, all in $V$.
After integration by parts, we obtain
\begin{subequations} \label{eqn:fem-3}
\begin{align}
\begin{split}
\int_\varOmega \psi \dx{u} \dx{x} + \int_\varOmega \psi u_x\left(u \dx{t} + \sum_j \varXi^j\circ \mathrm{d}W^j\right)\dx{x} \\ 
 - \int_\varOmega \psi_x F \dx{t} \dx{x}+\sum_j\int_\varOmega \psi G^j\circ \mathrm{d}W^j \dx{x}= & \ 0,
\end{split} \\
\int \phi F \dx{x} +\int_\varOmega \alpha^2 \phi_x F_{x}\dx{x} -\int_\varOmega \phi \left[u^2 + \frac{\alpha^2}{2}(u_x)^2\right]\dx{x} = & \ 0, \\
\int_\varOmega \zeta G^j \dx{x} + \int_\varOmega \alpha^2\zeta_x G^j_{x} - \int_\varOmega \zeta\left[ 2u\varXi^j_x + \alpha^2u_x\varXi^j_{xx}\right]\dx{x}= & \ 0.
\end{align}
\end{subequations}
To discretise in time, we use the implicit midpoint rule, defining
\begin{equation}
u^{(h)}=\frac{1}{2}\left(u^{(n+1)}+u^{(n)}\right),
\end{equation}
where $u^{(n)}$ is the value of $u$ at the $n$-th time level.
This is the natural choice for discretising a Stratonovich process and is unconditionally stable for Stratonvich stochastic differential equations (see \cite{abdulle2012high}).
This gives us our mixed finite element discretisation, in which we search for the simultaneous fields
\begin{equation}
\left(u^{(n+1)}, \Delta F^{(h)}, \Delta G^{(h)}\right)\in \left(V, V, V\right)
\end{equation}
that satisfy, for all $(\psi, \phi,\zeta)\in \left(V, V, V\right)$, the equations
\begin{subequations}\label{eqn:discretisation}
\begin{align}
& \int_\varOmega \psi\left(u^{(n+1)}-u^{(n)}\right) \dx{x} + \int_\varOmega \psi u_x^{(h)} \Delta v^{(h)} \dx{x} - \int_\varOmega \psi_x \Delta F^{(h)} \dx{x} +   \int_\varOmega \psi \Delta G^{(h)} \dx{x} = 0, \\
& \int_\varOmega \phi \Delta F^{(h)} \dx{x} + \alpha^2\int_\varOmega \phi_x \Delta F_x^{(h)} \dx{x} - \int_\varOmega \phi u^{(h)} u^{(h)} \Delta t \dx{x}
- \frac{\alpha^2}{2}\int_\varOmega \phi u_x^{(h)}u_x^{(h)}\Delta t \dx{x}=0, \\
& \int_\varOmega \zeta \Delta G^{(h)} \dx{x} + \alpha^2 \int_\varOmega \zeta_x \Delta G^{(h)}_x\dx{x} - \sum_j\left(\int_\varOmega 2\zeta u^{(h)}\varXi^j_x+\alpha^2 \zeta u_x^{(h)}\varXi^j_{xx}\right) \Delta W^j\dx{x}=0,
\end{align}
\end{subequations}
where $\Delta t$ is the time step and $\Delta W^j$ is a normally distributed number with zero mean and variance $\Delta t$.
The stochastic velocity $\Delta v^{(h)}$ is given by
\begin{equation}
\Delta v^{(h)} = u^{(h)}\Delta t + \sum_j \varXi^j \Delta W^j,
\end{equation}
We have also absorbed the time step $\Delta t$ and the increments $\Delta W$ into $F$ and $G^j$, writing these as $\Delta F$ and $\Delta G$ to emphasise that they contain increments in time.\\
\\
If the system is deterministic, there are no $\varXi^j$ functions that are non-zero and the equations reduce so that $\Delta G=0$ and $\Delta v^{(h)}=u^{(h)}\Delta t$.\\
\\
To implement this scheme, we used the Firedrake software described in \cite{rathgeber2017firedrake}.
\section{Convergence of the Discretisation to Peakon Solutions}\label{sec:numerical convergence}
\begin{figure}[!htb]
\centering
\includegraphics[width=\textwidth]{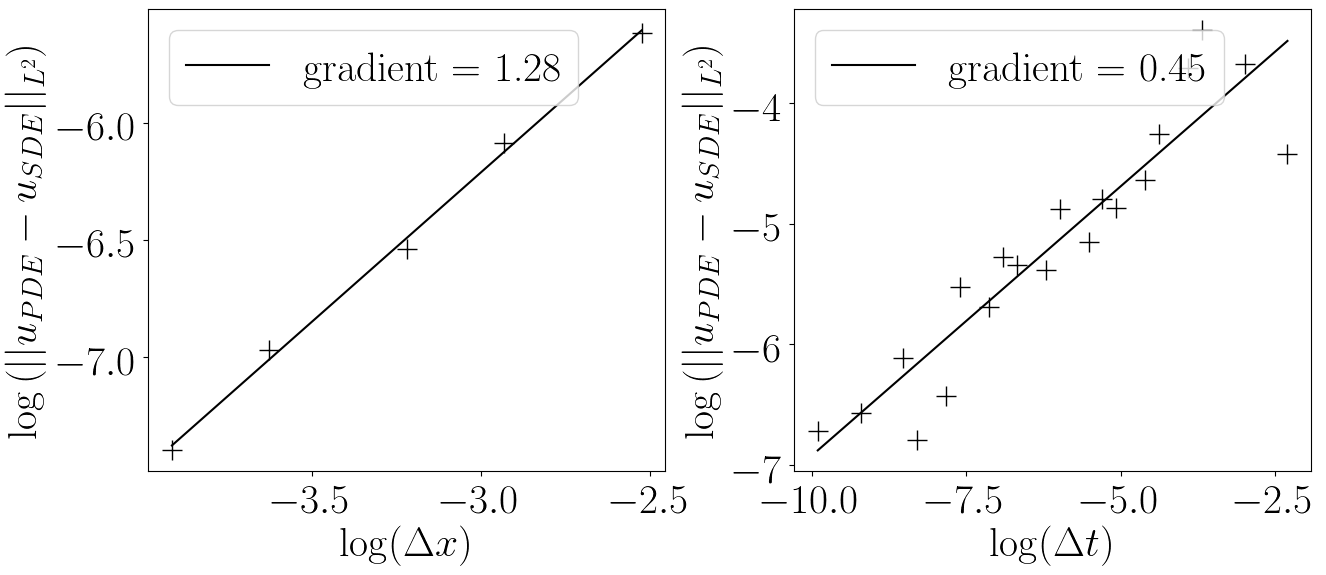}
\caption{Plots showing the convergence (left) as $\Delta x \to 0$ and (right) as $\Delta t\to 0$ of the discretisation from Section \ref{sec:discretisation} to the equations (\ref{eqn:diff peakons Ito}) that describe the evolution of peakon solutions.
The field $u_\mathrm{PDE}$ was computed by solving the discretisation from an initial condition resembling a peakon.
The equations (\ref{eqn:diff peakons Ito}) were solved using Milstein's method.
The resulting values of $p$ and $q$ were used to reconstruct the field $u_{\mathrm{SDE}}$.
We then computed the error between $u_\mathrm{PDE}$ and $u_{\mathrm{SDE}}$ using the $L^2$ norm.
Crosses mark the data points, with a solid best fit line overlaid. \\
\rule{\textwidth}{0.4pt}
} \label{fig:convergence}
\end{figure}
We showed in Section \ref{sec:weak peakons} that the peakon (\ref{eqn:periodic peakon}) satisfying (\ref{eqn:diff peakons Strat}) is a weak solution to (\ref{eqn:weak sch-hydro}).
In this section, we demonstrate that, when describing a single peakon, our discretisation of the partial differential equation (PDE) (\ref{eqn:sch hydro}) also converges to the stochastic differential equations (SDEs) (\ref{eqn:diff peakons Strat}) for the evolution of the canonical coordinates $p$ and $q$, in the limit that the grid spacing $\Delta x \to 0$.
We also demonstrate that the discretisation of the PDE has strong (pathwise) convergence to the peakon SDEs as the time step $\Delta t\to 0$. \\
\\
To do this, we first re-cast the equations (\ref{eqn:diff peakons Strat}) for a single peakon in It\^{o} form:
\begin{subequations}\label{eqn:diff peakons Ito}
\begin{align}
\dx{p} & = \frac{p}{2} \sum_j \left[\left(\varXi^j_x(q)\right)^2-\varXi^j(q)\varXi^j_{xx}(q)\right]\dx{t} -p \sum_j \varXi^j_x(q) \  \mathrm{d}W^j, \\
\dx{q} & = \frac{1}{2}\left[\frac{p(1+e^{-L/\alpha})}{1-e^{L/\alpha}}+\sum_j\varXi^j(q)\varXi^j_x(q)\right]\dx{t} + \sum_j \varXi^j(q)\dx{W}_j.
\end{align}
\end{subequations}
We solve both (\ref{eqn:discretisation}) and (\ref{eqn:diff peakons Ito}) in a periodic one-dimensional domain of length $L=40$, and take $\alpha = 1$.
For this section, for the stochastic basis functions we used only a single constant function:
\begin{equation}
\varXi^0 = 1.
\end{equation}
Our strategy is to establish a numerical solution of the SDEs (\ref{eqn:diff peakons Ito})
that is well-resolved temporally, by solving it using Milstein's method for discretising in time; see for instance \cite{jacobs2010stochastic} for more details of this scheme.
We used an initial condition of $p=1$ and $q=20$ and with $\Delta t=10^{-6}$.
Using the resulting $p$ and $q$, we can construct the velocity field $u$ corresponding to this peakon via (\ref{eqn:periodic peakon}).
This can be compared with the velocity field $u$ that has been computed by solving the stochastic PDE (\ref{eqn:sch hydro}) using the discretisation presented in Section \ref{sec:discretisation}.\\
\\
Our initial condition for the PDE, corresponding to $p=1$ and $q=20$ was
\begin{equation}
u = \left \lbrace \begin{matrix}
\frac{1}{2(1-e^{-L/\alpha})} \left[e^{(x-L/2)/\alpha} + e^{-L/\alpha}e^{-(x-L/2)/\alpha}\right], & x < L/2, \\
& \\
\frac{1}{2(1-e^{-L/\alpha})} \left[e^{-(x-L/2)/\alpha} + e^{-L/\alpha}e^{(x-L/2)/\alpha}\right], & x \geq L/2.
\end{matrix}\right.
\end{equation}
Taking different choices of $\Delta x$ and $\Delta t$, we computed the state $u_\mathrm{PDE}$ at $t=0.1$ for various spatial resolutions.
Then we measured the error in the $L^2$-norm between $u_\mathrm{PDE}$ and the corresponding $u_\mathrm{SDE}$ that was reconstructed from the $p$ and $q$ computed at $t=0.1$ from the numerical solution to the SDE.\\
\\
We did two separate investigations into the convergence of this error: firstly reducing $\Delta x$  with $\Delta t=10^{-6}$ kept constant, and secondly reducing $\Delta t$ with $\Delta x=0.002$ kept constant.
In the second case, as $\Delta t$ was changed it was important to ensure that the stochastic realisation still corresponded to that of `true' solution computed from the SDEs.
Let  $\curlybrackets{\Delta \widetilde{W}_1,\dots,\Delta \widetilde{W}_n}$ be the set of randomly generated increments used for the `true' solution, with $\Delta \widetilde{t}=10^{-6}$.
Then when taking a larger time step, such that $\Delta t = M \Delta \widetilde{t}$, we took
\begin{equation}
\Delta W_i = \frac{1}{\sqrt{M}}\sum_{j=1}^{M}\Delta \widetilde{W}_{(i-1)M+j}, 
\end{equation}
where the square root is necessary to keep the variance of $\Delta W$ equal to $\Delta t$.
A description of this approach can be found in \cite{jacobs2010stochastic}. \\
\\
The errors are plotted as functions of $\Delta x$ and $\Delta t$ in Figure \ref{fig:convergence}.
This demonstrates that when our discretisation presented in Section \ref{sec:discretisation} describes a peakon, it does indeed converge to the coupled SDEs (\ref{eqn:diff peakons Strat}) describing the evolution of a peakon.
Since the peakon profile does not have a continuous derivative, we do not expect to obtain a convergence rate approaching second-order as $\Delta x\to 0$.
Although a discretisation using the implicit midpoint rule should have first-order strong (pathwise) convergence as $\Delta t\to 0$, (see \cite{milstein2002numerical} or \cite{abdulle2012high}), again due to the lack of smoothness in the solution we expect the reduced convergence that we have observed.

\section{Numerical Investigations}\label{sec:results}
\subsection{Deterministic formation of peakons via wave breaking} \label{sec:results-deterministic-steep}
\begin{figure}[h!]
\centering
\includegraphics[width=\textwidth]{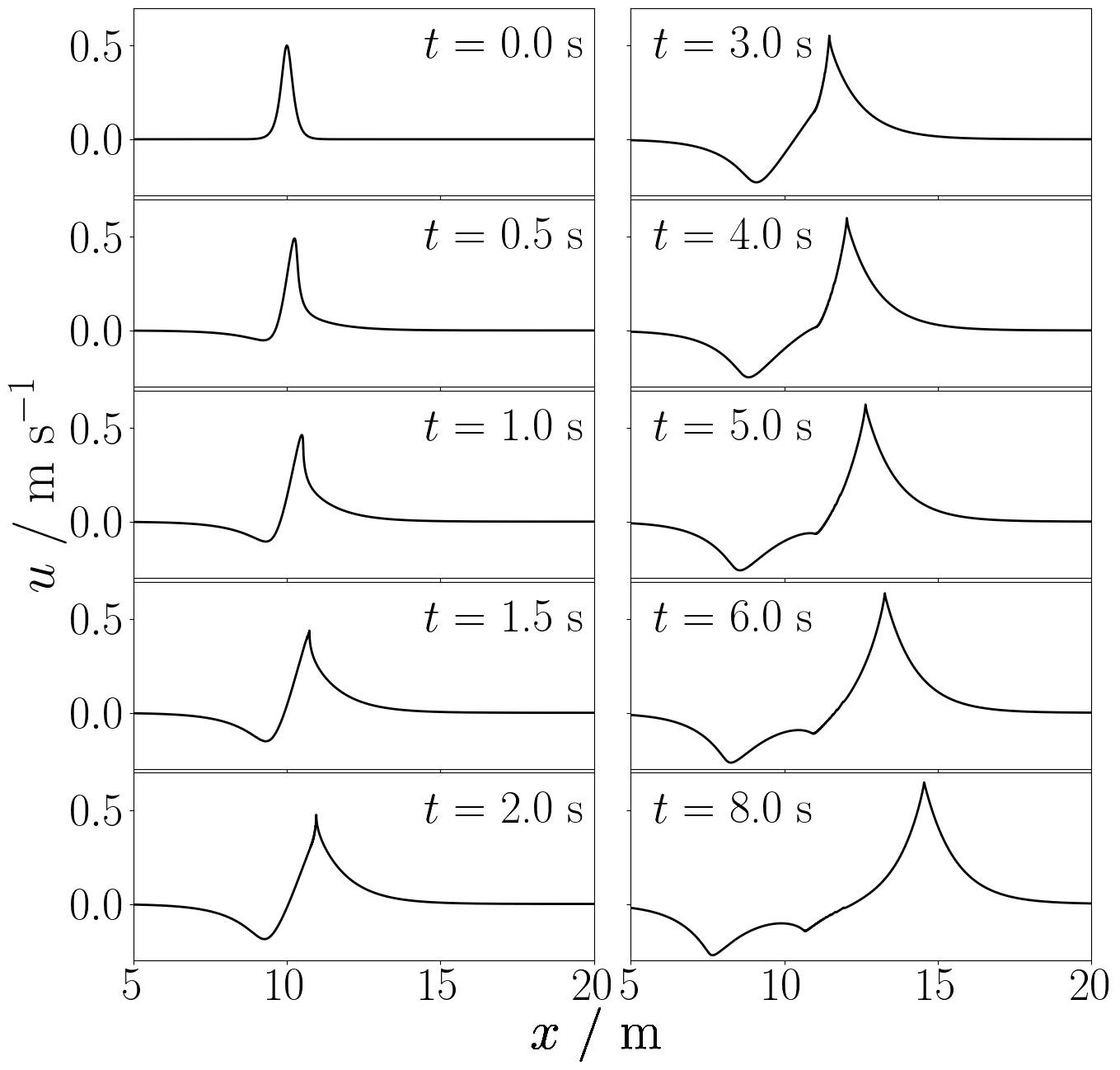}
\caption{The emergence of a peakon via wave breaking from the initial condition (\ref{eqn:steep initial condition}).
Each panel shows the velocity profile $u$ in a periodic domain (although only part of the domain is shown) at a different moment in time.
This initial condition has a momentum $m$ with negative parts and satisfies the assumptions of the steepening lemma of \cite{camassa1993integrable} and we see that the slope approaches vertical and a peakon forms through wave breaking.
In fact, we also see the formation of anti-peakons (peakons with negative velocity values) to the left of the primary peakon. \\
\rule{\textwidth}{0.4pt}}
\label{fig:steep demo}
\end{figure}
In this section we use the finite element discretisation to investigate the formation of peakons.
For all of Section \ref{sec:results} we take $\alpha=1$ and use a periodic domain of length $L=40$. \\
\\
First we will demonstrate the formation of a peakon via wave breaking in the deterministic case from a smooth, confined initial condition.
The steepening lemma of \cite{camassa1993integrable} considers profiles whose initial negative slope $s$ at the inflection point and the integral $H_1=\tfrac{1}{2}\int_\varOmega \left( u^2+\alpha^2u^2_x\right) \dx{x}$ satisfy $s < -\sqrt{2H_1/\alpha^3}$.
Here we choose an initial condition satisfying this:
\begin{equation}
u = \frac{1}{2} \ \mathrm{sech}\left(\frac{x-L/4}{L/240}\right), \label{eqn:steep initial condition}
\end{equation}
so that $s=-1.5$, $H_1\simeq 0.542$ and $\sqrt{2H_1/\alpha^3}\simeq 1.041$.
Computing $m=u-\alpha^2u_{xx}$ for this profile reveals $m$ to have a negative section to the right of a positive section, and as shown by \cite{mckean1998breakdown} and discussed by \cite{mckean2015breakdown} this is a condition for wave breaking.
The deterministic case is achieved in our discretisation by simply taking $\varXi=0$.
Figure \ref{fig:steep demo} shows the velocity profile found using the discretisation of Section \ref{sec:discretisation} with $\Delta x=0.002$ and $\Delta t=0.001$ at a range of times.
We see that the peak leans to the right and the negative slope steepens until it becomes vertical.
Then the wave breaks and a peakon emerges to propagate to the right.

\subsection{Deterministic formation of peakons without wave breaking} \label{sec:results-deterministic-shallow}
\begin{figure}[h!]
\centering
\includegraphics[width=\textwidth]{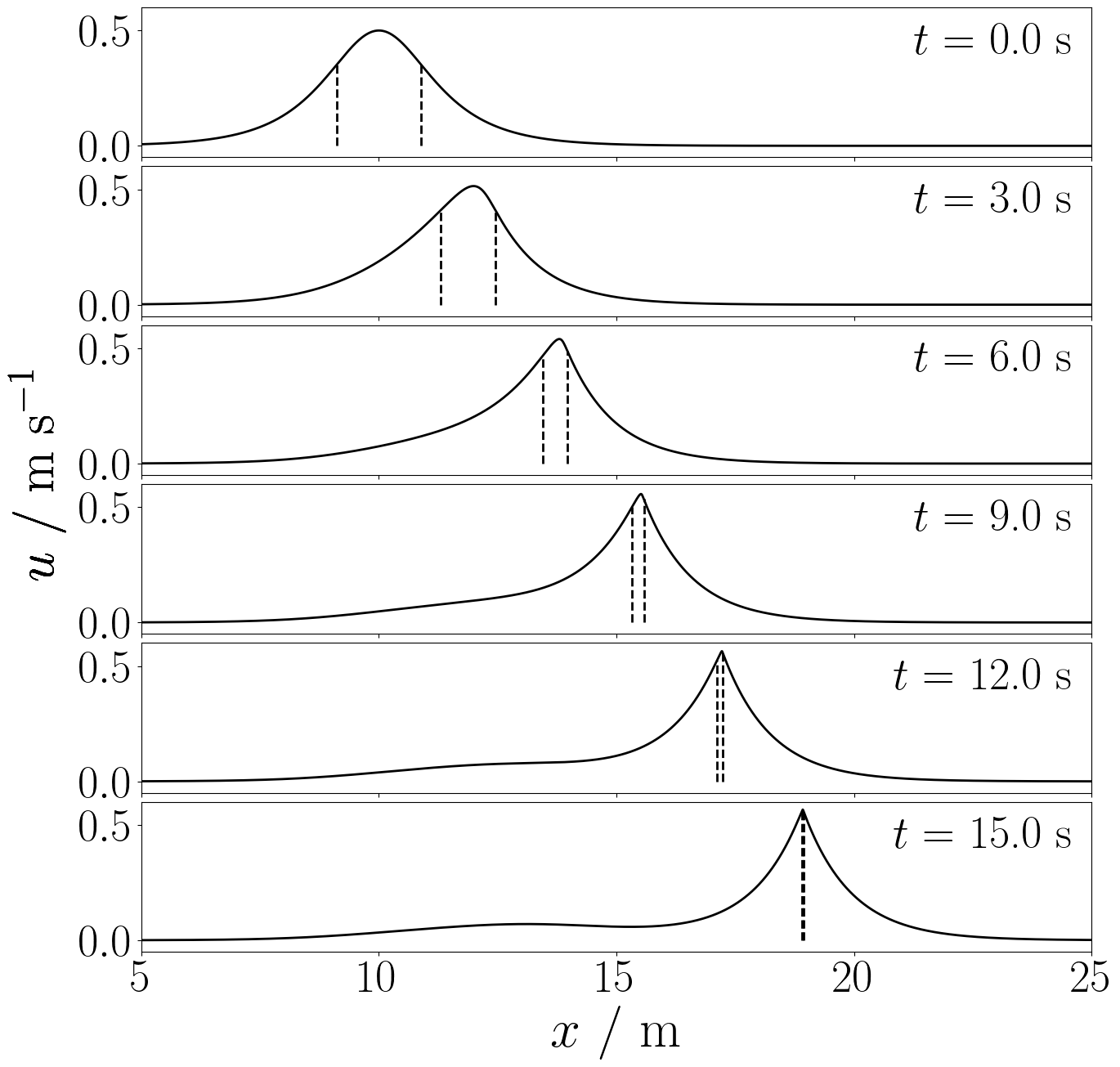}
\caption{The emergence of a peakon without wave breaking.
The initial condition used was (\ref{eqn:shallow initial condition}), whose slope is shallower than that assumed for the steepening lemma of \cite{camassa1993integrable} and whose momentum $m>0$.
As the profile moves rightwards, although we see it tilt the slope never becomes vertical.
Instead, the inflection points (marked with dashed lines) rise until they meet at the top of the profile.
At this moment, the peakon emerges. \\
\rule{\textwidth}{0.4pt}
}
\label{fig:shallow demo}
\end{figure}
Now we consider an initial condition with a shallower slope:
\begin{equation}
u = \frac{1}{2} \ \mathrm{sech}\left(\frac{x-L/4}{L/40}\right). \label{eqn:shallow initial condition}
\end{equation}
This does not satisfy the assumptions of the steepening lemma of \cite{camassa1993integrable}, as
$s=-0.25$, $H_1\simeq 0.333$ and $\sqrt{2H_1/\alpha^3}\simeq 0.816$.
In contrast to the steep initial condition in the previous section, computing $m=u-\alpha^2u_{xx}$ at each point gives a momentum that is entirely positive. This is a condition for no wave breaking discussed in \cite{mckean1998breakdown}. \\
\\
Figure \ref{fig:shallow demo} illustrates peakon formation from this initial condition.
It shows the velocity profile at a series of times, with the inflection points highlighted with dashed lines.
As the velocity is represented by continuous, piecewise linear function in our discretisation, the gradient $u_x$ of the field is a series of discontinuous piecewise constants.
The positions of the inflection points are then identified as the centre of the cells containing the maxima and minima in $u_x$.
As we see in Figure \ref{fig:shallow demo}, the inflection points rise up the profile until they eventually meet at the maximum -- which is the moment of peakon formation.
At no point does the slope turn vertical, and thus we see that peakon formation can occur without wave breaking.
\subsection{Diagnostics for peakon formation and wave breaking} \label{sec:results-diagnostics}
\begin{figure}[h!]
\centering
\includegraphics[width=\textwidth]{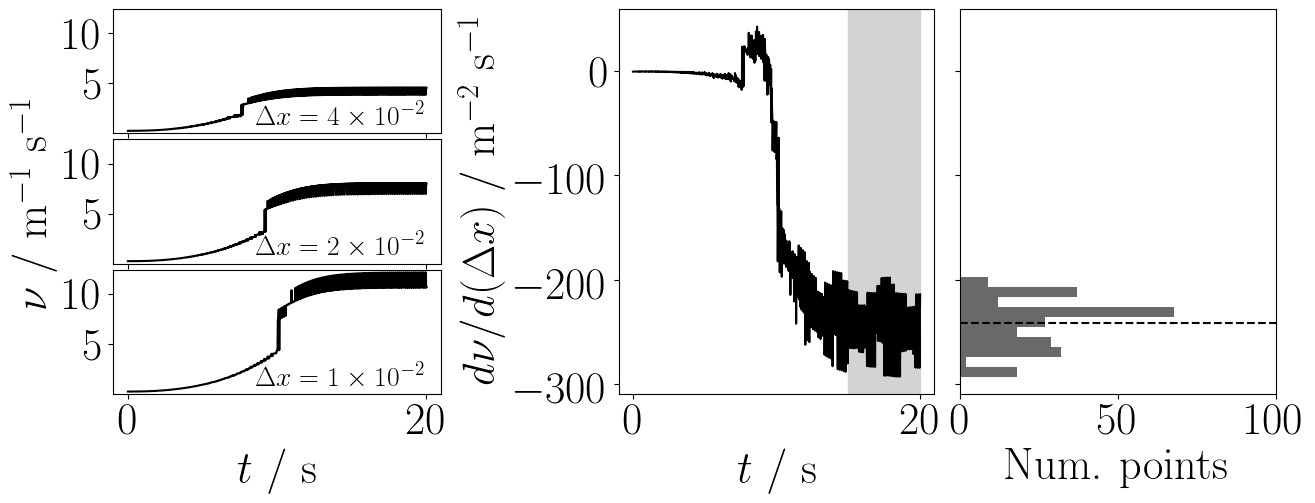}
\caption{An illustration of the process for determining whether a velocity profile is a peakon. This example uses the shallower initial condition (\ref{eqn:shallow initial condition}) of Section \ref{sec:results-deterministic-shallow}.
(Left) the diagnostic $\nu$ of (\ref{eqn:nu}) is plotted as a function of time for simulations at five different resolutions (three are shown).
We see that before the peakon forms, there is agreement in the values of $\nu$ which has converged to some finite value.
Then there is a clear jump in the value of $\nu$, and as the resolution is refined this jump gets larger, and also at a later point in time.
After the peakon has formed, there is no longer convergence in $\nu$ as the resolution is refined.
The values oscillate as the peakon moves through a cell, giving the resulting thick band.
(Centre) at each point in time, a linear curve is fitted to the values of $\nu$ as a function of $\Delta x$.
The central panel plots the gradient of this curve as a function of time, which we denote by $\mathrm{d}\nu/\mathrm{d}(\Delta x)$.
We see that initially this value is close to $0$ as the value of $\nu$ has converged.
After the peakon has formed, there is a dramatic drop in $\mathrm{d}\nu/\mathrm{d}(\Delta x)$, showing that as the grid spacing is reduced the value of $\nu$ gets increasingly higher.
To determine whether the peakon has formed, we inspect the values of $\mathrm{d}\nu/\mathrm{d}(\Delta x)$ from $t\in[15,20]$ s.
These times are shaded in grey on the central figure.
(Right) the values of $\mathrm{d}\nu/\mathrm{d}(\Delta x)$ at these times are plotted as a histogram, with the mean value shown with the dashed line, giving the value $\Pi$ of equation (\ref{eqn:Pi}). \\
\rule{\textwidth}{0.4pt}}
\label{fig:deterministic nu shallow}
\end{figure}
\begin{figure}[h!]
\centering
\begin{subfigure}{\textwidth}
\centering
\includegraphics[width=\textwidth]{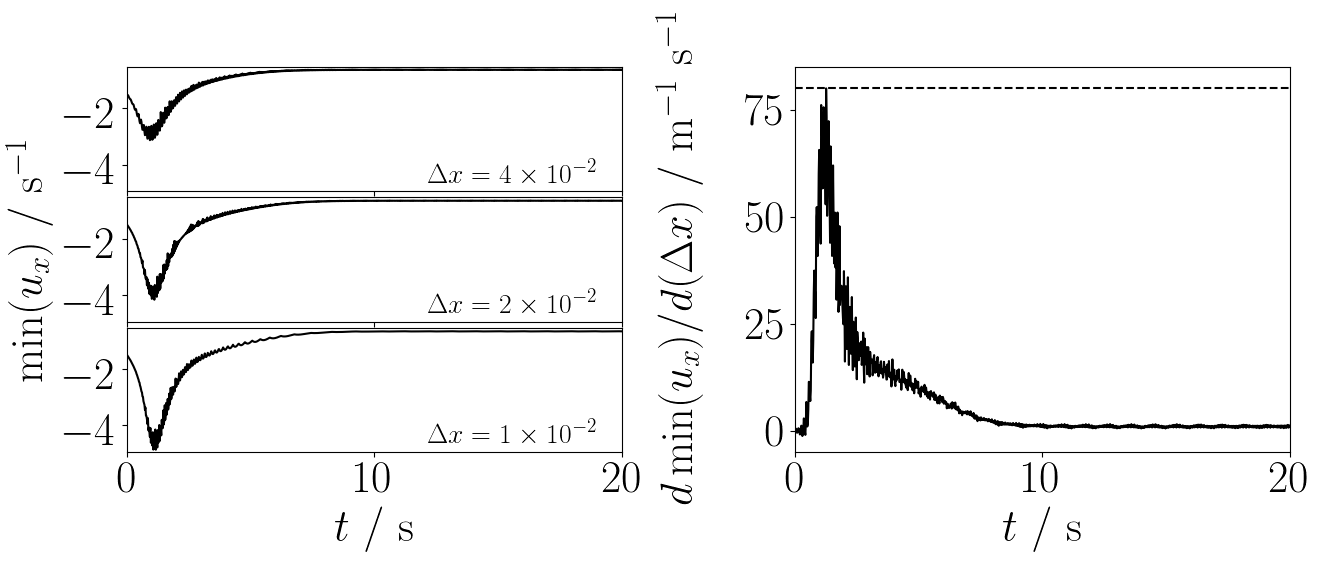}
\end{subfigure}
\begin{subfigure}{\textwidth}
\centering
\includegraphics[width=\textwidth]{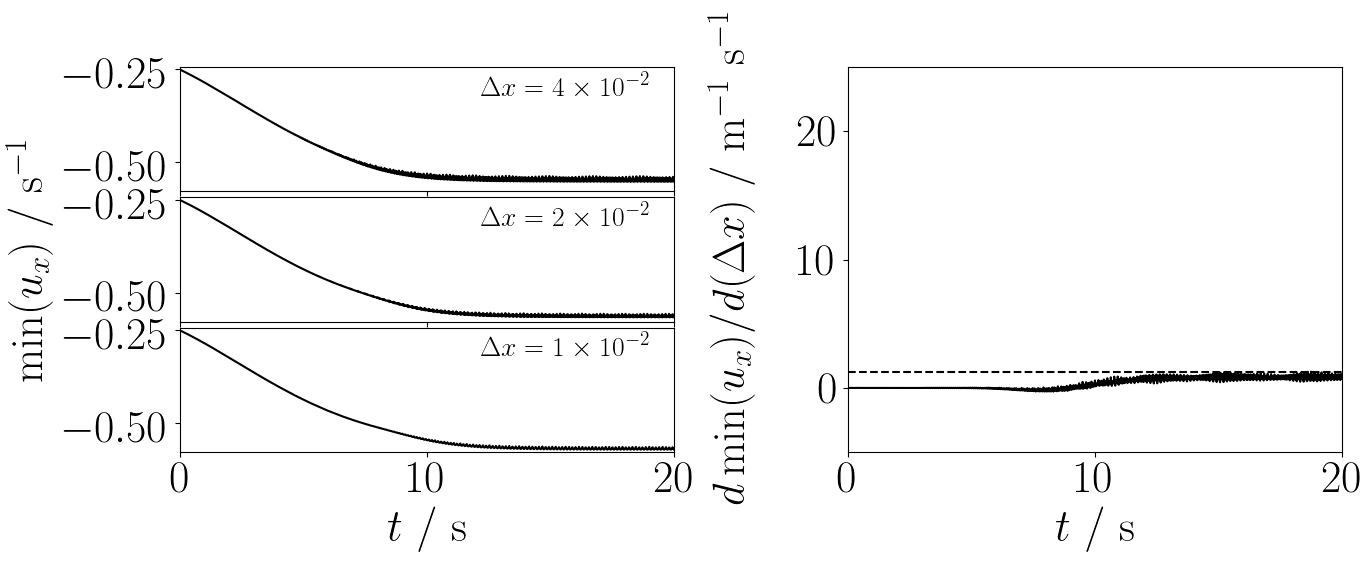}
\end{subfigure}
\caption{An illustration for the process of determining whether wave breaking has taken place. The top figures show the process for the steep initial condition of Section \ref{sec:results-deterministic-steep}, while the bottom figures show the process for the shallow profile of Section \ref{sec:results-deterministic-shallow}.
(Left) we plot $\min(u_x)$ as a function of time for five different resolutions (three are shown for each case).
For the steeper initial condition (top left), we see that roughly between $t=1$ s and $t=3$ s there is a dramatic drop in $\min(u_x)$, which appears steeper as the resolution is refined, signifying that the slope of $u$ is turning vertical.
(Right) for each point in time we fit a linear curve to $\min(u_x)$ as a function of $\Delta x$ and plot the gradients of these curves, denoting this by $\mathrm{d}\min(u_x)/\mathrm{d}(\Delta x)$. As expected we see a strong positive here, representing the wave breaking. The dashed black line shows the maximum value.
(Bottom) in contrast, for the shallower initial profile we see agreement in $\min(u_x)$ for the different resolutions and values of $\mathrm{d}\min(u_x)/\mathrm{d}(\Delta x)$ that are close to zero. This demonstrates a lack of wave breaking. \\
\rule{\textwidth}{0.4pt}}
\label{fig:deterministic mindu}
\end{figure}
Before we investigate the formation of peakons in the stochastic case, we first define two diagnostics which allow us to identify whether a velocity profile is a peakon, and whether wave breaking has taken place.\\
\\
Our first diagnostic, which we denote by $\nu$, describes the jump in $u_x$ around the peak and is given by
\begin{equation} \label{eqn:nu}
\nu := \frac{\max_x(u_x)-\min_x(u_x)}{x_{\min}-x_{\max}},
\end{equation}
where $x_{\max}$ and $x_{\min}$ are the locations of the maximum and minimum in $u_x$.
For a peakon, as the grid spacing $\Delta x\to 0$, $\nu\to\infty$, whereas it will tend to some finite value when not representing a peakon.\\
\\
Figure \ref{fig:deterministic nu shallow} shows the evolution with time of $\nu$ for three values of $\Delta x$, using the shallow set-up outlined in Section \ref{sec:results-deterministic-shallow}.
The time step used was $\Delta t=0.0005$, although when this was repeated with smaller values of $\Delta t$ there was no significant difference in the values.
We see that as the slope on the leading edge of the velocity profile steepens, the value of $\nu$ rises, but there is agreement as the resolution is refined and $\nu$ tends to some finite value.
Then the values of $\nu$ jump sharply, and as expected it appears that as $\Delta x\to 0$, $\nu\to\infty$.
The subsequent band of oscillating values corresponds to the changing value of $(x_{\max}-x_{\min})$ as the peakon moves through a cell.
To detect the presence of a peakon, at each point in time we inspect $\mathrm{d}\nu/\mathrm{d}(\Delta x)$ by fitting a linear curve to $\nu$ as a function of $\Delta x$ and inspecting the gradient.
If this is strongly negative it represents the presence of a peakon.
This is also displayed in Figure \ref{fig:deterministic nu shallow}.
As we will run a large quantity of simulations with different stochastic realisations, it is useful to return a single value indicating the presence of a peakon.
For this, we look at the average value in $\mathrm{d}\nu/\mathrm{d}(\Delta x)$ for times close to the end of the simulation, which for brevity we will call:
\begin{equation} \label{eqn:Pi}
\Pi := \left\langle \frac{\mathrm{d}\nu}{\mathrm{d}(\Delta x)} \right\rangle,
\end{equation}
with the average denoted by $\left\langle\cdot\right\rangle$ taken over $t\in[15, 20]$ s.
Thus the presence of a peakon is signified by a large negative value of $\Pi$.
The lack of a peakon would be represented by a value of $\Pi$ close to zero.
\\
\\
To detect wave breaking, we look at $\min_x(u_x)$ as a function of time.
If wave breaking is taking place and the slope turns vertical, then as $\Delta x\to 0$ we will see $\min_x(u_x)\to -\infty$.
Again for each point in time we fit a linear curve to $\min(u_x)$ as a function of $\Delta x$ and inspect the gradient, with strongly positive values indicating wave breaking.
Figure \ref{fig:deterministic mindu} shows this for the set-ups of Sections \ref{sec:results-deterministic-steep} and \ref{sec:results-deterministic-shallow}, clearly showing wave breaking occurring with the steeper initial condition of equation (\ref{eqn:steep initial condition}) but not that of equation (\ref{eqn:shallow initial condition}).
As before, it is helpful to define a single value to determine whether wave breaking has taken place.
We therefore look at 
\begin{equation}
\omega:=\max_t\left[\frac{\mathrm{d}\min_x(u_x)}{\mathrm{d}(\Delta x)} \right],
\end{equation}
with large positive values of $\omega$ representing wave breaking and values of $\omega$ close to zero showing a lack of wave breaking.

\subsection{Stochastic formation of peakons} \label{sec:results-stochastic} % under constant $\varXi$}
\begin{figure}
\centering
\includegraphics[width=\textwidth]{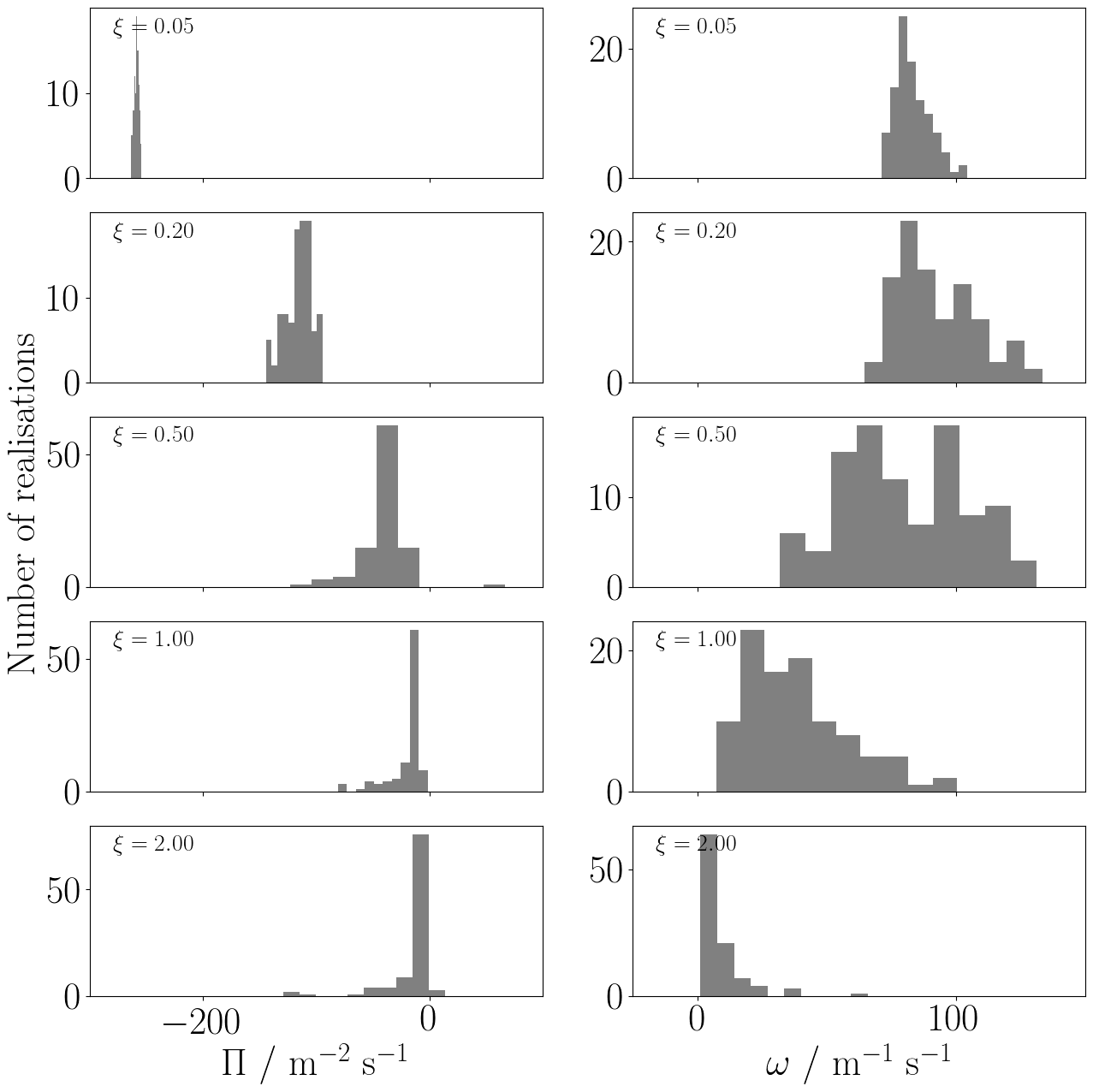}
\caption{Histograms of the peakon formation diagnostic $\Pi$ and the wave breaking diagnostic $\omega$ for stochastic simulations from the steep initial condition of Section \ref{sec:results-deterministic-steep}.
The different rows show different strengths of stochastic perturbation, each with 100 realisations of the noise.
Negative values of $\Pi$ indicate peakon formation, while positive values of $\omega$ represent wave breaking having occurred.
These results signal that peakon formation and wave breaking can still occur in the stochastic Camassa-Holm equation. \\
\rule{\textwidth}{0.4pt}}
\label{fig:steep stochastic hist}
\end{figure}
\begin{figure}
\centering
\includegraphics[width=\textwidth]{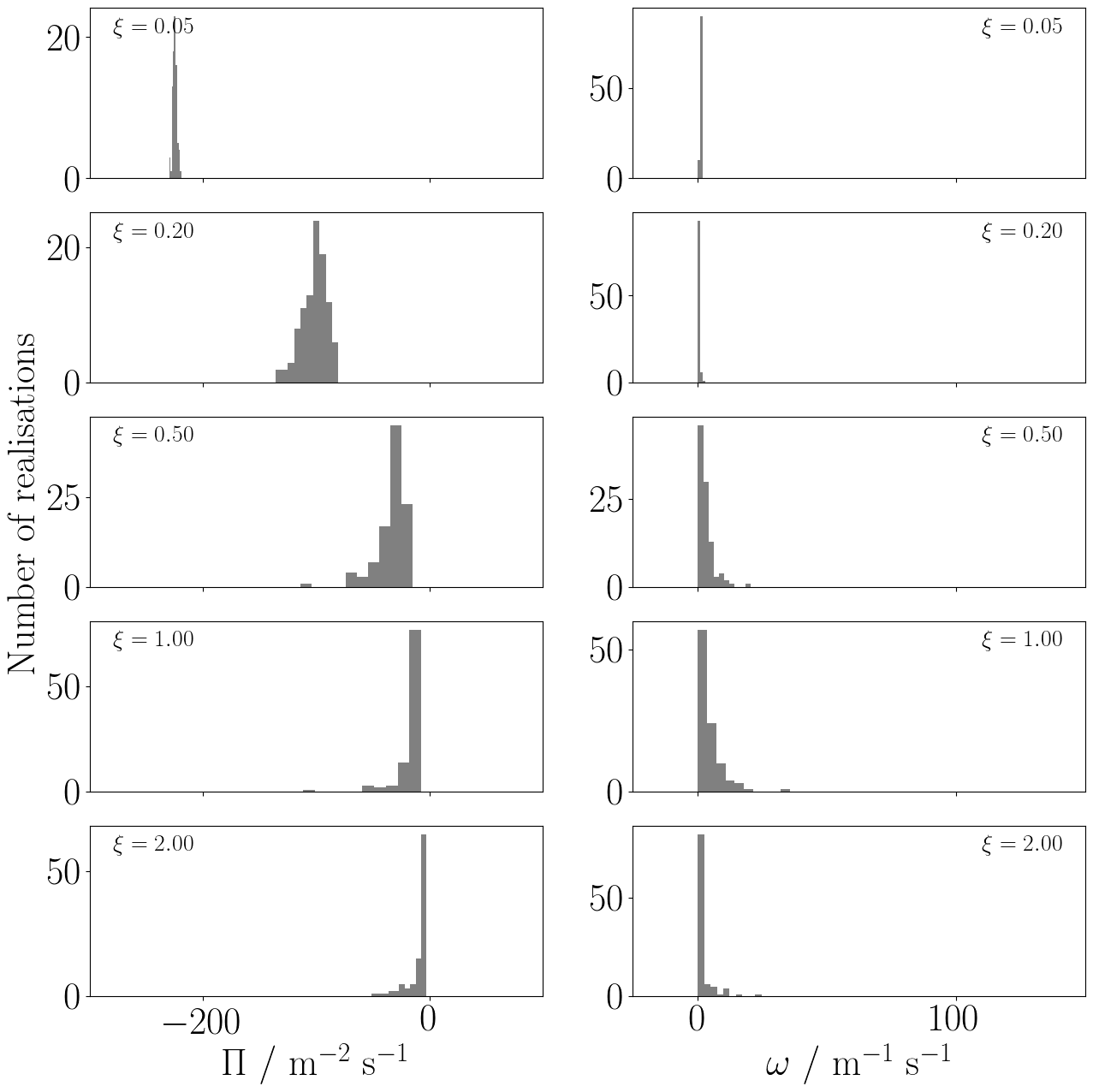}
\caption{Histograms of the peakon formation diagnostic $\Pi$ and the wave breaking diagnostic $\omega$ for stochastic simulations from the steep initial condition of Section \ref{sec:results-deterministic-steep}.
The different rows show different strengths of stochastic perturbation, each with 100 realisations of the noise.
Negative values of $\Pi$ indicate peakon formation, while positive values of $\omega$ represent wave breaking having occurred.
These results support the claim that peakon formation does occur in the stochastic Camassa-Holm equation, and that it can occur without wave breaking taking place. \\
\rule{\textwidth}{0.4pt}}
\label{fig:shallow stochastic hist}
\end{figure}
In \cite{crisan2018wave}, the authors investigated the formation of peakons with the case that $\varXi(x)=\xi$, a constant.
They did this by considering a stochastic form of the steepening lemma of \cite{camassa1993integrable}, applied to initial conditions satisfying the same condition of $s<-\sqrt{2H_1/\alpha^3}$.
The authors found that wave breaking did occur with a positive probability, but it was unclear whether this probability was necessarily unity.
Although this was interpreted as the probability of peakon formation, the results from Sections \ref{sec:results-deterministic-steep} to \ref{sec:results-diagnostics} indicate that instead this was the probability of the wave breaking process occurring.\\
\\
In this section, we use the discretisation of Section \ref{sec:discretisation} to investigate both wave breaking and the formation of peakons in the SCH equation.
First we aim to determine qualitatively whether wave breaking and peakon formation may still occur following the introduction of stochastic transport.
Finding that they do, we may ask whether the probabilities of these events is merely non-zero or in fact unity.
To answer these questions, we use the diagnostics $\Pi$ and $\omega$ and the methodologies that were laid out in Section \ref{sec:results-diagnostics}.
Large negative values of $\Pi$ indicate the formation of a peakon, while large positive values of $\omega$ indicate that wave breaking has taken place.\\
\\
For both the steep and shallow initial profiles laid out in Sections \ref{sec:results-deterministic-steep} and \ref{sec:results-deterministic-shallow}, we now add a stochastic component to the transport with $\varXi=\xi$ (a constant), for a range of $\xi$ values.
Simulations were all computed up to $t=20$ s.
To investigate the probabilities of peakon formation and wave breaking, we solved the SCH equation under 100 different realisations of the noise, for each value of $\xi$ and both the steep and shallow initial profiles.
In each case we computed values of the diagnostics $\Pi$ and $\omega$, which are plotted as histograms for each $\xi$ in Figures \ref{fig:steep stochastic hist} and \ref{fig:shallow stochastic hist}.\\
\\
For each realisation of the noise we performed the simulation at five different resolutions, in each case dividing our domain uniformly into $1000$, $1500$, $2000$, $2500$ and $3000$ cells.
The values of $\nu$ and $\min_x(u_x)$ at each of these resolutions were used to determine the diagnostics $\Pi$ and $\omega$.
We used a time step of $\Delta t=0.0005$ s. As shown in Section \ref{sec:numerical convergence}, for a stochastic peakon the convergence as $\Delta t\to 0$ is slow.
Thus we found that refining $\Delta t$ further did improve the results, but limits on computational resources then restricted the number of realisations of noise that it was feasible to use.\\
\\
We look first at Figure \ref{fig:steep stochastic hist}, which uses the steep initial condition from Section \ref{sec:results-deterministic-steep}.
The values of $\Pi$ are negative in every single case, indicating that peakons can indeed form.
Two positive values of $\Pi$ were recorded, though further inspection of the fields from these simulations suggested that these were still peakons.
Similarly the values of $\omega$ are positive in every case, particularly strongly when the noise is weak, suggesting that wave breaking can occur from the initial condition of equation (\ref{eqn:steep initial condition}).\\
\\
Figure \ref{fig:shallow stochastic hist} plots the diagnostics $\Pi$ and $\omega$ for the shallow initial condition of Section \ref{sec:results-deterministic-shallow}.
As with the steeper initial condition, the results suggest that peakon formation does take place.
However we did not see a clear example of wave breaking, although it is possible that for large values of $\xi$ that wave breaking may take place.\\
\\
Both Figures \ref{fig:steep stochastic hist} and \ref{fig:shallow stochastic hist} show that as the strength of the noise is increased, the signals of both peakon formation and wave breaking get weaker, as the diagnostics of $\nu$ and $\min_x(u_x)$ become more noisy.
As the simulations had not converged with respect to $\Delta t$, we did see improvements in these signals by reducing $\Delta t$ further, but as mentioned above this results in a trade-off with the number of realisations of the noise.
However, even for the strongest values of $\xi$ we do see evidence of peakon formation.\\
\\
It is also worth noting that our $\Pi$ diagnostic could show false positive results for velocities profiles that are not peakons, but whose inflection points are separated by less than a distance of the order of the smallest $\Delta x$ that we used.
However, we have demonstrated in Section \ref{sec:weak peakons} that peakons are weak solutions to the SCH equation and in Section \ref{sec:numerical convergence} that our discretisation will converge to such solutions.\\
\\
The conclusions that we draw must come with caveats.
There are infinitely many initial conditions that we have not considered, with infinitely many choices of $\varXi$ and infinitely many realisations of the noise.
It may be that peakons might not form under larger values of $\xi$, although the physical motivation of including stochastic transport to represent unresolved processes suggests that values of $\xi$ should be small.
Despite this, we did not see any examples in which a peakon did not appear to form.
In this section we have therefore seen that under the addition of stochastic perturbations both peakon formation and wave breaking can still occur.

\section{Conclusions}
We have considered the Camassa-Holm equation of \cite{camassa1993integrable} with the addition of the stochastic transport of \cite{holm2015variational}.
Taking the advective form of this stochastic Camassa-Holm equation, we showed that peakons are weak solutions to this equation in Section \ref{sec:weak peakons}.
Inspired by looking for weak solutions, we presented a finite element discretisation to the advective form of the equation in Section \ref{sec:discretisation} and demonstrated that it converges to the peakon solutions in Section \ref{sec:numerical convergence}.
The strength of using such a finite element discretisation is that it can be used to describe solutions of the equation through the process of peakon formation. \\
\\
The simulations presented in Section \ref{sec:results} demonstrate the formation of peakons through two mechanisms.
Firstly through wave breaking, as described by the steepening lemma of \cite{camassa1993integrable} when the slope at the inflection point to the right of the initial peak is sufficiently negative.
For shallower initial profiles, peakons can still form without wave breaking taking place. 
Instead, the inflection points rise up the profile without the slope ever turning vertical and the peakon emerges as the inflection points reach the summit.\\
\\
We then investigated these processes with stochastic perturbations using a constant $\varXi$, testing the findings of \cite{crisan2018wave}.
Our results suggest that peakon formation does still occur, and that wave breaking can still occur when the initial slope is steep.
We did not see any examples when peakons did not form in these situations.\\
\\
A clear avenue for further work is to investigate whether spatially-varying $\varXi$ may prevent the formation of peakons.
The (impossibly unlikely) choice of $\sum_j\varXi^j\circ \mathrm{d}W^j=-u \dx{t}$ in equation (\ref{eqn:SCH-eqn}) yields $\dx{u}=0$, suggesting that this may be possible.
We have also not discussed the effect that the noise has upon the time of the peakon formation.
\section*{Acknowledgements and Funding}
All three authors were supported by EPSRC grant EP/N023781/1.
TMB was also supported by the EPSRC Mathematics of Planet Earth Centre for Doctoral Training at Imperial College London and the University of Reading, with grant number EP/L016613/1.
The authors would like to thank the three anonymous reviewers, whose comments and suggestions helped improve the manuscript.
\bibliography{peakon_formation}

\begin{thebibliography}{10}

\bibitem{bertini1994stochastic}
L.~Bertini, N.~Cancrini, and G.~Jona-Lasinio, ``The stochastic {B}urgers
  equation,'' {\em Communications in Mathematical Physics}, vol.~165, no.~2,
  pp.~211--232, 1994.

\bibitem{alonso2019burgers}
D.~Alonso-Or{\'a}n, A.~B. de~Le{\'o}n, and S.~Takao, ``The {B}urgers' equation
  with stochastic transport: shock formation, local and global existence of
  smooth solutions,'' {\em Nonlinear Differential Equations and Applications
  NoDEA}, vol.~26, no.~6, p.~57, 2019.

\bibitem{flandoli2010well}
F.~Flandoli, M.~Gubinelli, and E.~Priola, ``Well-posedness of the transport
  equation by stochastic perturbation,'' {\em Inventiones mathematicae},
  vol.~180, no.~1, pp.~1--53, 2010.

\bibitem{flandoli2011random}
F.~Flandoli, {\em Random Perturbation of {PDE}s and Fluid Dynamic Models:
  {\'E}cole d'{\'e}t{\'e} de Probabilit{\'e}s de Saint-Flour XL--2010},
  vol.~2015.
\newblock Springer Science \& Business Media, 2011.

\bibitem{camassa1993integrable}
R.~Camassa and D.~D. Holm, ``An integrable shallow water equation with peaked
  solitons,'' {\em Physical Review Letters}, vol.~71, no.~11, p.~1661, 1993.

\bibitem{holm2016variational}
D.~D. Holm and T.~M. Tyranowski, ``Variational principles for stochastic
  soliton dynamics,'' {\em Proceedings of the Royal Society A: Mathematical,
  Physical and Engineering Sciences}, vol.~472, no.~2187, 2016.

\bibitem{crisan2018wave}
D.~Crisan and D.~D. Holm, ``Wave breaking for the stochastic {C}amassa--{H}olm
  equation,'' {\em Physica D: Nonlinear Phenomena}, vol.~376, pp.~138--143,
  2018.

\bibitem{holm2015variational}
D.~D. Holm, ``Variational principles for stochastic fluid dynamics,'' {\em
  Proceedings of the Royal Society A: Mathematical, Physical and Engineering
  Sciences}, vol.~471, no.~2176, 2015.

\bibitem{bendall2018statistical}
T.~M. Bendall and C.~J. Cotter, ``Statistical properties of an enstrophy
  conserving finite element discretisation for the stochastic quasi-geostrophic
  equation,'' {\em Geophysical \& Astrophysical Fluid Dynamics}, pp.~1--14,
  2018.

\bibitem{cotter2019numerically}
C.~Cotter, D.~Crisan, D.~D. Holm, W.~Pan, and I.~Shevchenko, ``Numerically
  modeling stochastic {L}ie transport in fluid dynamics,'' {\em Multiscale
  Modeling \& Simulation}, vol.~17, no.~1, pp.~192--232, 2019.

\bibitem{constantin1998global}
A.~Constantin and J.~Escher, ``Global existence and blow-up for a shallow water
  equation,'' {\em Annali della Scuola Normale Superiore di Pisa-Classe di
  Scienze}, vol.~26, no.~2, pp.~303--328, 1998.

\bibitem{dullin2001integrable}
H.~R. Dullin, G.~A. Gottwald, and D.~D. Holm, ``An integrable shallow water
  equation with linear and nonlinear dispersion,'' {\em Physical Review
  Letters}, vol.~87, no.~19, p.~194501, 2001.

\bibitem{dullin2003camassa}
H.~R. Dullin, G.~A. Gottwald, and D.~D. Holm, ``{C}amassa--{H}olm,
  {K}orteweg--de {V}ries-5 and other asymptotically equivalent equations for
  shallow water waves,'' {\em Fluid Dynamics Research}, vol.~33, no.~1-2,
  p.~73, 2003.

\bibitem{dullin2004asymptotically}
H.~Dullin, G.~Gottwald, and D.~Holm, ``On asymptotically equivalent shallow
  water wave equations,'' {\em Physica D: Nonlinear Phenomena}, vol.~190,
  no.~1-2, pp.~1--14, 2004.

\bibitem{constantin2001scattering}
A.~Constantin, ``On the scattering problem for the {C}amassa-{H}olm equation,''
  {\em Proceedings of the Royal Society of London. Series A: Mathematical,
  Physical and Engineering Sciences}, vol.~457, no.~2008, pp.~953--970, 2001.

\bibitem{fuchssteiner1981symplectic}
B.~Fuchssteiner and A.~S. Fokas, ``Symplectic structures, their {B}{\"a}cklund
  transformations and hereditary symmetries,'' {\em Physica D: Nonlinear
  Phenomena}, vol.~4, no.~1, pp.~47--66, 1981.

\bibitem{fuchssteiner1996some}
B.~Fuchssteiner, ``Some tricks from the symmetry-toolbox for nonlinear
  equations: generalizations of the {C}amassa-{H}olm equation,'' {\em Physica
  D: Nonlinear Phenomena}, vol.~95, no.~3-4, pp.~229--243, 1996.

\bibitem{chen2012well}
Y.~Chen, G.~Hongjun, and G.~Boling, ``Well--posedness for stochastic
  {C}amassa--{H}olm equation,'' {\em Journal of Differential Equations},
  vol.~253, no.~8, pp.~2353--2379, 2012.

\bibitem{lv2020dependence}
G.~Lv, J.~Wei, and G.-a. Zou, ``The dependence on initial data of stochastic
  {C}amassa--{H}olm equation,'' {\em Applied Mathematics Letters}, p.~106472,
  2020.

\bibitem{constantin1998wave}
A.~Constantin and J.~Escher, ``Wave breaking for nonlinear nonlocal shallow
  water equations,'' {\em Acta Mathematica}, vol.~181, no.~2, pp.~229--243,
  1998.

\bibitem{constantin2000global}
A.~Constantin and L.~Molinet, ``Global weak solutions for a shallow water
  equation,'' {\em Communications in Mathematical Physics}, vol.~211, no.~1,
  pp.~45--61, 2000.

\bibitem{holm2005momentum}
D.~D. Holm and J.~E. Marsden, ``Momentum maps and measure-valued solutions
  (peakons, filaments, and sheets) for the {EPD}iff equation,'' in {\em The
  breadth of symplectic and Poisson geometry}, pp.~203--235, Springer, 2005.

\bibitem{bressan2007global}
A.~Bressan and A.~Constantin, ``Global conservative solutions of the
  {C}amassa--{H}olm equation,'' {\em Archive for Rational Mechanics and
  Analysis}, vol.~183, no.~2, pp.~215--239, 2007.

\bibitem{holden2007global}
H.~Holden and X.~Raynaud, ``Global conservative solutions of the
  {C}amassa--{H}olm equation - a {L}agrangian point of view,'' {\em
  Communications in Partial Differential Equations}, vol.~32, no.~10,
  pp.~1511--1549, 2007.

\bibitem{holm2018stochastic}
D.~D. Holm and T.~M. Tyranowski, ``Stochastic discrete hamiltonian variational
  integrators,'' {\em BIT Numerical Mathematics}, vol.~58, no.~4,
  pp.~1009--1048, 2018.

\bibitem{holm2018new}
D.~D. Holm and T.~M. Tyranowski, ``New variational and multisymplectic
  formulations of the euler--poincar{\'e} equation on the virasoro--bott group
  using the inverse map,'' {\em Proceedings of the Royal Society A:
  Mathematical, Physical and Engineering Sciences}, vol.~474, no.~2213,
  p.~20180052, 2018.

\bibitem{mckean1998breakdown}
H.~McKean, ``Breakdown of a shallow water equation,'' {\em Asian Journal of
  Mathematics}, vol.~2, no.~4, pp.~867--874, 1998.

\bibitem{constantin1999shallow}
A.~Constantin and H.~P. McKean, ``A shallow water equation on the circle,''
  {\em Communications on Pure and Applied Mathematics: A Journal Issued by the
  Courant Institute of Mathematical Sciences}, vol.~52, no.~8, pp.~949--982,
  1999.

\bibitem{xu2008local}
Y.~Xu and C.-W. Shu, ``A local discontinuous {G}alerkin method for the
  {C}amassa--{H}olm equation,'' {\em SIAM Journal on Numerical Analysis},
  vol.~46, no.~4, pp.~1998--2021, 2008.

\bibitem{holden2008numerical}
H.~Holden and X.~Raynaud, ``A numerical scheme based on multipeakons for
  conservative solutions of the {C}amassa--{H}olm equation,'' in {\em
  Hyperbolic problems: theory, numerics, applications}, pp.~873--881, Springer,
  2008.

\bibitem{holden2006convergence}
H.~Holden and X.~Raynaud, ``Convergence of a finite difference scheme for the
  {C}amassa--{H}olm equation,'' {\em SIAM journal on numerical analysis},
  vol.~44, no.~4, pp.~1655--1680, 2006.

\bibitem{artebrant2006numerical}
R.~Artebrant and H.~J. Schroll, ``Numerical simulation of {C}amassa--{H}olm
  peakons by adaptive upwinding,'' {\em Applied numerical mathematics},
  vol.~56, no.~5, pp.~695--711, 2006.

\bibitem{chertock2012convergence}
A.~Chertock, J.-G. Liu, and T.~Pendleton, ``Convergence of a particle method
  and global weak solutions of a family of evolutionary {PDE}s,'' {\em SIAM
  Journal on Numerical Analysis}, vol.~50, no.~1, pp.~1--21, 2012.

\bibitem{liu2016invariant}
H.~Liu and Y.~Xing, ``An invariant preserving discontinuous {G}alerkin method
  for the {C}amassa--{H}olm equation,'' {\em SIAM Journal on Scientific
  Computing}, vol.~38, no.~4, pp.~A1919--A1934, 2016.

\bibitem{li2014high}
M.~Li and A.~Chen, ``High order central discontinuous {G}alerkin-finite element
  methods for the {C}amassa--{H}olm equation,'' {\em Applied Mathematics and
  Computation}, vol.~227, pp.~237--245, 2014.

\bibitem{matsuo2010hamiltonian}
T.~Matsuo, ``A {H}amiltonian-conserving {G}alerkin scheme for the
  {C}amassa--{H}olm equation,'' {\em Journal of computational and applied
  mathematics}, vol.~234, no.~4, pp.~1258--1266, 2010.

\bibitem{antonopoulos2019error}
D.~Antonopoulos, V.~A. Dougalis, and D.~Mitsotakis, ``Error estimates for
  {G}alerkin finite element methods for the {C}amassa--{H}olm equation,'' {\em
  Numerische Mathematik}, vol.~142, no.~4, pp.~833--862, 2019.

\bibitem{abdulle2012high}
A.~Abdulle, D.~Cohen, G.~Vilmart, and K.~C. Zygalakis, ``High weak order
  methods for stochastic differential equations based on modified equations,''
  {\em SIAM Journal on Scientific Computing}, vol.~34, no.~3, pp.~A1800--A1823,
  2012.

\bibitem{rathgeber2017firedrake}
F.~Rathgeber, D.~A. Ham, L.~Mitchell, M.~Lange, F.~Luporini, A.~T. McRae, G.-T.
  Bercea, G.~R. Markall, and P.~H. Kelly, ``Firedrake: automating the finite
  element method by composing abstractions,'' {\em ACM Transactions on
  Mathematical Software (TOMS)}, vol.~43, no.~3, p.~24, 2017.

\bibitem{jacobs2010stochastic}
K.~Jacobs, {\em Stochastic processes for physicists: understanding noisy
  systems}.
\newblock Cambridge University Press, 2010.

\bibitem{milstein2002numerical}
G.~N. Milstein, Y.~M. Repin, and M.~V. Tretyakov, ``Numerical methods for
  stochastic systems preserving symplectic structure,'' {\em SIAM Journal on
  Numerical Analysis}, vol.~40, no.~4, pp.~1583--1604, 2002.

\bibitem{mckean2015breakdown}
H.~P. McKean, ``Breakdown of the {C}amassa-{H}olm equation,'' in {\em Henry P.
  McKean Jr. Selecta}, pp.~189--193, Springer, 2015.

\end{thebibliography}
\bibliographystyle{ieeetr}

\end{document}